\documentclass[portuguese,a4paper,twoside,10pt]{article}
\usepackage[english]{babel}
\usepackage[T1]{fontenc}
\usepackage[latin1]{inputenc}
\usepackage{amssymb,amsmath,epsfig,amsfonts,euscript}
\usepackage{dsfont}
\usepackage[usenames,dvipsnames]{color}

\newtheorem{teo}{Theorem}[section]
\newtheorem{pro}[teo]{Proposition}

\newcounter{note}[section]
\newcounter{example}[section]

\newcommand{\dem}{{\bf Proof. }}
\newcommand{\fdem}{ $\square$}

\def\N{{\rm I\kern-.20em N}}
\def\R{{\rm I\kern-.20em R}}
\def\indi{{1\kern-.20em\rm I}}

\newcommand{\titulo}[1]{\begin{center}\mbox{} \\ \noindent \textit{\textbf{\Large #1}}\\\vspace{0.5cm}\end{center}}

\newcommand{\autor}[1]{\noindent \textbf{#1}}

\newcommand{\afil}[1]{{\small \noindent \textit{#1}}}

\renewcommand{\abstract}[1]{{\small \noindent \textbf{Abstract:} #1\\}}

\newcommand{\keywords}[1]{{\small \noindent \textbf{Keywords:} #1\\}}

\newcommand{\biblio}[1]{{\small \begin{thebibliography}{99} #1 \end{thebibliography}}}

\begin{document}
% ---------------------------------------------------------------------
% T\'{\i}tulo do artigo.
% ---------------------------------------------------------------------
\titulo{Generalized Logistic Models and its orthant tail dependence}

\autor{Helena Ferreira, Luisa Pereira}

\afil{Department of Mathematics, University
of Beira Interior, Covilh\~a, Portugal}\\

\abstract{The Multivariate Extreme Value distributions have shown their usefulness in environmental studies, financial and insurance mathematics. The Logistic or Gumbel-Hougaard distribution is one of the oldest multivariate extreme value models and it has been extended to asymmetric models. In this paper we introduce generalized logistic multivariate distributions. Our tools are mixtures of copulas and stable mixing variables, extending approaches in Tawn (1990), Joe and Hu (1996) and Foug\`eres et al. (2009). The parametric family of multivariate extreme value distributions considered presents a flexible dependence structure and we compute for it the multivariate tail dependence coefficients considered in Li (2009).}


%%%%%%%%%%%%%%%%%%%%%%%%%%%%%%%%%%%%%%%%%%%%%%%%%%%%%%%
\section{The model}

Let $\mathcal{L}\left(Z|W\right)$ denotes the conditional distribution of a random variable or vector $Z$ given another random variable or vector $W$. For the vectors ${\bf X}_j=\left(X_{j,1},\ldots,X_{j,d}\right)$, $j=1,...,q$, and ${\bf S}=\left(S_1,\ldots,S_q\right)$, defined on the same probability space, we shall assume that:\\

(a) $\mathcal{L}\left(\left({\bf X}_1,...,{\bf X}_q\right)|{\bf S}\right)=\prod_{j=1}^q\mathcal{L}\left({\bf X}_j|{\bf S}\right)$,\\

(b)$\mathcal{L}\left({\bf X}_j|{\bf S}\right)=\mathcal{L}\left({\bf X}_j|S_j\right)$\\

(c) $P\left(\bigcap_{i=1}^dX_{ji}\leq x_i|S_j\right)=C_j\left(e^{-\left(\frac{x_1}{\beta_{j1}}\right)^{-1/\alpha_j}S_j},...,e^{-\left(\frac{x_d}{\beta_{jd}}\right)^{-1/\alpha_j}S_j}\right)$, $x_j>0$, $j=1,...,d$, where $C_j$'s are max-stable copulas and $\{\beta_{ji}, j=1,...,q, i=1,...,d\}$ are non-negative constants such that $\sum_{j=1}^q\beta_{ji}=1$, $ i=1,...,d$,\\

(d) $E\left(e^{-tS_j}\right)=e^{-t^{\alpha_j}}$, $t\geq 0$, $j=1,...,q$, where $\alpha_j$'s are constants in $(0,1]$\\

and\\

(e) $\mathcal{L}\left({\bf S}\right)=\prod_{j=1}^q\mathcal{L}\left(S_j\right)$.\\

Thus every $X_{ji}$ is a scale mixture with mixing variable $\beta_{ji}S_j^{\alpha_j}$ and ${\bf X}_j$, $j=1,...,q$, are conditionally independent given ${\bf S}$.\\
Scale mixtures have been studied and used in a variety of applications ( Marshall and Olkin (1988, \cite{marshall+olkin}), Joe and Hu (1996, \cite{joe+hu}) and Foug\`eres et al. (2009, \cite{fougeres+nolan+rootzen}), Li (2009, \cite{li})).\\

We shall consider here a componentwise maxima model from the ${\bf X}_j$'s. From this model we derive a new family of copulas and analize its orthant tail dependence by computing the multivariate tail dependence coefficients considered in Li (2009, \cite{li})). Finally we apply the results to the particular case of $C_j$ being the copula arising from the distribution of the variables in a M4 process (Smith and Weissman, 1996, \cite{smith+weissman}).

%*************************************************************************************************************************
 \begin{pro}\label{p1.1}
If the random vectors  ${\bf X}_j$, $j=1,...,q$, and ${\bf S}$ satisfy the conditions (a)-(e) then ${\bf Y}=\left(Y_1,\ldots,Y_d\right)$ defined by $Y_i=\displaystyle\bigvee_{j=1}^q X_{ji}$, $i=1,...,d$, has multivariate extreme value distribution with unit Frechet margins and copula
\begin{eqnarray}\label{modelo}
C_{{\bf Y}}\left(u_1,...,u_d\right)=\exp\left\{-\sum_{j=1}^q\left(-\ln C_j\left(e^{-\left(-\beta_{j1}\ln u_1\right)^{1/\alpha_j}},...,e^{-\left(-\beta_{jd}\ln u_d\right)^{1/\alpha_j}}\right)\right)^{\alpha_j}\right\}.
\end{eqnarray}

\end{pro}
 
 \dem To obtain $C_{{\bf Y}}$ we just apply the conditional independence of  the ${\bf X}_j$'s follwoed by the max-stability of $C_j$'s and the $\alpha_j$-stability of each $S_j$, as follows:
 
 $$P\left(\bigcap_{i=1}^d\left\{Y_i\leq x_i\right\}\right)=\int P\left(\bigcap_{j=1}^q\bigcap_{i=1}^d\left\{X_{ji}\leq x_i\right\}|S=s\right)d{\bf S}\left(s_1,...,s_q\right)=$$
 
 $$\int \prod_{j=1}^qC_j\left(e^{-\left(\frac{x_1}{\beta_{j1}}\right)^{-1/\alpha_j}s_j},...,e^{-\left(\frac{x_d}{\beta_{jd}}\right)^{-1/\alpha_j}s_j}\right)d{\bf S}\left(s_1,...,s_q\right)=$$
 
 $$\prod_{j=1}^q\exp\left\{-\left(-\ln C_j\left(e^{-\left(\frac{x_1}{\beta_{j1}}\right)^{-1/\alpha_j}},...,e^{-\left(\frac{x_d}{\beta_{jd}}\right)^{-1/\alpha_j}}\right)\right)^{\alpha_j}\right\}.$$
 
 For each $j$ and $i$, $X_{ji}$ is a positive $\alpha_j$-stable size mixture of a Frechet distribution with location $\beta_{ji}$, scale $\beta_{ji}\alpha_j$ and shape parameter $\alpha_j$ and has itself Frechet distribution with same location and the same right end point, but scale $\beta_{ji}$ and shape parameter $1$. Since $\sum_{j=1}^q\beta_{ji}=1$, $ i=1,...,d$, each $Y_i$ has unit Frechet distribution. The max-stability of $C_{{\bf Y}}$ follows from its expression and the max-stability of the  $C_j$'s.\fdem
 \vspace{0.4cm}

%The family of copulas $C_{{\bf Y}}$ is closed under margins, that is, for positive constants $c_i$, $i=1,...,d$, the variable $\bigvee_{i=1}^dc_iY_i$ is a random variable in the same scale family as the $Y_i$'s.

We now discuss some particular cases of (\ref{modelo}) that has been explored.\\

%**********************
{\bf (I)} If $q=1$ then $\beta_{1i}=1$,  $i=1,...d$, and 

$$C_{{\bf Y}}\left(u_1,...,u_d\right)=\exp\left\{-\left(-\ln C_1\left(e^{-\left(-\ln u_1\right)^{1/\alpha_1}},...,e^{-\left(-\ln u_d\right)^{1/\alpha_1}}\right)\right)^{\alpha_1}\right\}$$

is a generalisation of the  Archimedean copula (Joe, 1997, \cite{joe}), which for the particular case of the product copula $C_1=\Pi$ leads to the Gumbel-Hougaard or logistic copula. The dependence properties of the special case of $C_1(u_1,...,u_d)=\displaystyle\prod_{1\leq s<t\leq d}C_{\{s,t\}}(u_s^{p_s},u_t^{p_t})\prod_{i=1}^du_i^{p_i\nu_i}$, where $C_{\{s,t\}}$, $1\leq s<t\leq d$, are bivariate copulas and $(d-1)p_i+p_i\nu_i=1$, $i=1,...,d$, were analysed in Joe and Hu (1996, \cite{joe+hu}).\\

{\bf (II)} If $C_j=\Pi$, $j=1,...,q$, then 

$$\prod_{j=1}^q\exp\left\{-\left(-\ln C_j\left(e^{-\left(\frac{x_1}{\beta_{j1}}\right)^{-1/\alpha_j}},...,e^{-\left(\frac{x_d}{\beta_{jd}}\right)^{-1/\alpha_j}}\right)\right)^{\alpha_j}\right\}=
\exp\left\{-\sum_{j=1}^q\left(\sum_{i=1}^d  \left(\frac{x_i}{\beta_{ji}}\right)^{-1/\alpha_j}    \right)^{\alpha_j}\right\},$$

 which leads to an asymmetric logistic copula 
 
\begin{eqnarray}\label{modelo1}
C_{{\bf Y}}\left(u_1,...,u_d\right)=\exp\left\{-\sum_{j=1}^q\left(\sum_{i=1}^d  \left(-\beta_{ji}\ln u_i\right)^{-1/\alpha_j}    \right)^{\alpha_j}\right\}.
\end{eqnarray}

In (\ref{modelo1}), if we take $\alpha_j=\alpha$, $j=1,...,q\leq +\infty$, we find an analogous mixture of extreme value distributions to those considered in Fougeres et al (2009, \cite{fougeres+nolan+rootzen}) by departing just from a random vector ${\bf X}=\left(X_1,\ldots,X_d\right)$ satisfying $\mathcal{L}\left({\bf X}_j|{\bf S}\right)=\prod_{i=1}^d\mathcal{L}\left(X_i|{\bf S}\right)$ and 
$P\left(X_i\leq x|{\bf S}\right)=\exp\{-\left(\sum_{j=1}^qc_{ji}S_j\right)\left(1+\gamma_i\frac{x-\mu_i}{\sigma_i}\right)^{-1/\gamma_i}\}$, $i=1,...,d$. That this, in this different approach, conditionally on ${\bf S}$, the vector  ${\bf X}$ has independent margins and each margin is a power mixture of an extreme value distribution with mixing variable $\sum_{j=1}^qc_{ji}S_j$, where the $c_{ji}$ are non-negartive constants.\\

{\bf (III)} Assume now, in (\ref{modelo1}), that each $j$ corresponds to an element $A$ of the set $\mathcal{S}$, the class of all nonempty subsets of $D=\{1,...,d\}$. If $\beta_{Ai}=0$ for each $i\not\in A$ then the copula  (\ref{modelo1}) becomes 

\begin{eqnarray}\label{modelo2}
C_{{\bf Y}}\left(u_1,...,u_d\right)=\exp\left\{-\sum_{A\subset \mathcal{S}}\left(\sum_{i\in A}^d  \left(-\beta_{Ai}\ln u_i\right)^{-1/\alpha_A}    \right)^{\alpha_A}\right\},
\end{eqnarray}
with $\sum_{A\subset \mathcal{S}}\beta_{Ai}=1$, $ i=1,...,d$. This is the asymmetric logistic model  considered in Tawn (1990, \cite{tawn}), by following a different probabilistic approach.
More generally, by applying the same interpretation of the constants $\beta_{ji}$ in (\ref{modelo}), we obtain 

\begin{eqnarray}\label{modelo3}
C_{{\bf Y}}\left(u_1,...,u_d\right)=\exp\left\{\sum_{A\subset \mathcal{S}}\left(-\ln C_A\left(e^{-\left(-\beta_{Ai_1(A)}\ln u_1\right)^{1/\alpha_A}},...,e^{-\left(-\beta_{Ai_s(A)}\ln u_d\right)^{1/\alpha_A}}\right)\right)^{\alpha_A}\right\},
\end{eqnarray}
where $C_A$'s are copulas with different dimensions and we denote by $\left(i_1(A),...,i_s(A)\right)$ the sub-vector of $(1,...,d)$ corresponding to indices in $A$. In particular, if we begin  with one copula $C_j=C$, $j=1,...,q$, then $C_A$, $A\subset \mathcal{S}$, are all the sub-copulas of $C$.\\

{\bf (IV)} Finally, let us suppose that $\beta_{ji}=\beta_j$, $ i=1,...,d$, in (\ref{modelo}). Then 

\begin{eqnarray}\label{modelo4}
C_{{\bf Y}}\left(u_1,...,u_d\right)=\prod_{j=1}^q\exp\left\{-\left(\ln C_j\left(e^{-\left(-\ln u_1\right)^{1/\alpha_j}},...,e^{-\left(-\ln u_d\right)^{1/\alpha_j}}\right)\right)^{\alpha_j}\beta_j\right\},
\end{eqnarray}
with $\sum_{j=1}^q\beta_j=1$, that is, $C_{{\bf Y}}$ is a geometric mean of mixtures of powers of multivariate extreme value distributions. The particular case of the weighted geometric mean $C_{{\bf Y}}(u_1,u_2)=(u_1\wedge u_2)^{\beta_1} (u_1u_2)^{1-\beta_1}$ is due to Cuadras and Aug\'e (1981, \cite{cuadras+auge}).
%=====================================================================================
\section{Orthant tail dependence}

For a random vector ${\bf Y}=\left(Y_1,\ldots,Y_d\right)$ with continuos margins $F_1,...,F_d$ and copula $C$, let the bivariate (upper) tail dependence  parameters defined by 

\begin{eqnarray}%\label{}
\lambda_{\{s,t\}}^{\left({\bf Y}\right)}\equiv \lambda_{\{s,t\}}^{\left(C\right)} =\lim_{u\uparrow 1} P\left(F_s\left(Y_s\right)>u|F_{t}\left(Y_{t}\right)>u\right),\,\,1\leq s<t\leq d.
\end{eqnarray}
The tail dependence is a copula based measure and it holds 
\begin{eqnarray} \label{1}
 \lambda_{\{s,t\}}^{\left(C\right)} =2-\lim_{u\uparrow 1} \frac{\ln C_{\{s,t\}}\left(u,u\right)}{\ln u},
\end{eqnarray}

where $C_{\{s,t\}}$ is the copula of the sub-vector $\left(Y_s,Y_{t}\right)$ ( Joe (1997, \cite{ joe}), Nelsen (1999, \cite{nelsen})).

To characterise the relative strength of extremal dependence with respect to a particular subset of random variables of ${\bf Y}$ one can use conditional orthant tail probabilities  of ${\bf Y}$ given that the components with indices in the subset $J$ are extreme. The tail dependence of bivariate copulas can be extended as done in Schmid and Schmidt (2007) (\cite{schmid/schmidt}) and Li (2009) (\cite{li}).

For $\emptyset \not=J\subset D=\{1,...,d\}$, let 

\begin{eqnarray}\label{lambda}
\lambda_{J}^{\left({\bf Y}\right)}\equiv \lambda_{J}^{\left(C\right)} =\lim_{u\uparrow 1} P\left(\displaystyle\bigcap_{j\notin J}\left\{F_j\left(Y_j\right)>u\right\}|
\displaystyle\bigcap_{j\in J}\left\{F_{j}\left(Y_{j}\right)>u\right\}\right).
\end{eqnarray}

If for some  $\emptyset \not=J\subset \{1,...,d\}$ the parameter $\lambda_{J}^{\left(C\right)}$ exists and is positive then we say that ${\bf Y}$ is (upper) orthant tail dependent. 

We have $\lambda_{J}^{\left(C\right)}=\frac{\lambda_{\{s\}}^{\left(C\right)}}{\lambda_{\{s\}}^{\left(C_{J}\right)}}$,  if $\lambda_{\{s\}}^{\left(C_{J}\right)}\not= 0$ and the relation  (\ref{1}) between the tail dependence parameter and the bivariate copula can also be generalized by

\begin{eqnarray}\label{lamdacopula}
 \lambda_{J}^{\left(C\right)} =\lim_{u\uparrow 1} \frac{\displaystyle\sum_{\emptyset \not=A\subset D}(-1)^{|A|-1}\ln C_A\left({\bf u}_A\right)}
 {\displaystyle\sum_{\emptyset \not=A\subset J}(-1)^{|A|-1}\ln C_A\left({\bf u}_A\right)},
\end{eqnarray}

where $C_A$ denotes the sub-copula of $C$ corresponding to margins with indices in $A$ and and ${\bf u}_A$ the $|A|$-dimensional vector $(u,...,u)$. By applying this relation and the max-stability of the copulas $C_j$, we get the following result.

 \begin{pro}\label{p2.1} 
For a copula $C$ defined by (\ref{modelo}), it holds

\begin{itemize}
\item[(a)]  
\begin{eqnarray}\label{lambda1}
\lambda_J^{\left(C\right)}=\frac
{\displaystyle\sum_{j=1}^q\displaystyle\sum_{\emptyset\not=A\subset D}(-1)^{|A|-1}\left(-\ln C_{j,A}\left(e^{-\beta_{j1}^{1/\alpha_j}},...,e^{-\beta_{jd}^{1/\alpha_j}}\right)_A\right)^{\alpha_j}}
{\displaystyle\sum_{j=1}^q\displaystyle\sum_{\emptyset\not=A\subset J}(-1)^{|A|-1}\left(-\ln C_{j,A}\left(e^{-\beta_{j1}^{1/\alpha_j}},...,e^{-\beta_{jd}^{1/\alpha_j}}\right)_A\right)^{\alpha_j}},
\end{eqnarray}
where $C_{j,A}$ denotes the sub-copula of $C_j$ corresponding to the margins with indices in $A$.\\

\item[(b)]  If $C_j=\Pi$, for each $j=1,...,q$, then
\begin{eqnarray}\label{lambda2}
\lambda_J^{\left(C\right)}=\frac
{\displaystyle\sum_{j=1}^q\displaystyle\sum_{\emptyset\not=A\subset D}(-1)^{|A|-1}\left(\displaystyle\sum_{i\in A }\beta_{ji}^{1/\alpha_j}\right)^{\alpha_j}}
{\displaystyle\sum_{j=1}^q\displaystyle\sum_{\emptyset\not=A\subset J}(-1)^{|A|-1}\left(\displaystyle\sum_{i\in A }\beta_{ji}^{1/\alpha_j}\right)^{\alpha_j}}.
\end{eqnarray}
\end{itemize}
\end{pro}
 
 \vspace {0.2cm}
 
 The tail dependence result in (\ref{lambda1}) depends on the mixing variables through the parameters $\alpha_j$, even for the case of $q=1$, that is the global dependence added by the mixing variables doesn't  vanish in extremes of maxima. This contrast with the result in Li (2009, \cite{li}), where the scale mixture of MEV distributions $(RX_1,...,RX_d)$ is considered with the mixing variable $R$ satisfying $\frac{E(e^{-ctR})}{E(e^{-tR})}\rightarrow c^{-\alpha}$,  as $t$ tends to $\infty$, and $c\geq 1$, $\alpha>0$. In this case the upper tail dependence parameters are exactly the same as the parameters of the MEV distribution without mixing.\\

 We remark that, for $\beta_{ji}=\beta_j$, $ i=1,...,d$, the numerator in   (\ref{lambda1}) is, for each $A\subset D$, 
 $$ \lambda_{\{s\}}^{\left(C_A\right)}=\sum _{j=1}^q\beta_j\lambda_{\{s\}}^{\left(C_{j,A}\right)}$$
 that is, the tail dependence parameter $ \lambda_{\{s\}}^{\left(C_A\right)}$ is a linear convex combination of the corresponding tail dependence parameters for the sub-copulas $C_{j,A}$ of $C_j$, $j=1,...,q$.\\

 The result in (\ref{lambda2}) leads to 
$$ \lambda_{\{s,t\}}^{\left(C\right)}=2-\sum_{j=1}^q\left(\beta_{js}^{1/\alpha_j}+\beta_{jt}^{1/\alpha_j}\right)^{\alpha_j},$$

extending the the known result 

\begin{eqnarray}\label{lambda3}
\lambda_{\{s,t\}}^{\left(C\right)}=2-2^\alpha,
\end{eqnarray}

 corresponding to $q=1$ (Joe (1997, \cite{ joe}), Nelsen (1999, \cite{nelsen})). The result in (\ref{lambda1}) enables to extend the equation (\ref{lambda3})  for other copulae $C_1$ than the product copula as 
 
 \begin{eqnarray}\label{lambda4}
\lambda_{\{s,t\}}^{\left(C\right)}=2-(2-\lambda_{\{s,t\}}^{\left(C_1\right)})^\alpha.
\end{eqnarray}

\section{Example}

We will suppose that $C_j=C$, $j=1,...,d$, with $$C(u_1,...,u_d)=
\prod_{l=1}^{\infty}\prod_{k=-\infty}^{\infty}\left(\displaystyle\bigwedge_{i=1}^d u_i^{a_{lki}}\right),\,\,\,u_j \in [0,1],\,\,j=1,...,d,$$
%\wedge ...\wedge u_d^{\alpha_{lkd}}\right),\,\,\, 
%$u_j \in [0,1]$, $j=1,...,d$, 

where $\{a_{lkj}, l\geq 1, -\infty<k<\infty,1\leq
j\leq d\},$ are nonnegative constants satisfying
$$\sum_{l=1}^\infty
\sum_{k=-\infty}^\infty a_{lkj}=1\quad {\rm for}\
j=1,\ldots,d.$$ That copula arises from the common distribution of the variables of an M4 process (Smith and Weissman (1996, \cite{smith+weissman})). 

Then the copula in (\ref{modelo}) becomes 

\begin{eqnarray}\label{exemplo1}
C_{{\bf Y}}\left(u_1,...,u_d\right)=\exp\left\{-\sum_{j=1}^q\left(
\sum_{l=1}^{\infty} \sum_{k=-\infty}^{\infty} \displaystyle\bigvee_{i=1}^d\left(-\beta_{ji}a_{lki}^{\alpha_j}\ln u_i\right)^{1/\alpha_j}
\right)^{\alpha_j}\right\}.
\end{eqnarray}

By applying the result in Proposition 2.1. (a), we obtain for the numerator in (\ref{lambda1})

$$\lambda_{\{s\}}^{\left(C_{{\bf Y}}\right)}=\displaystyle\sum_{j=1}^q\displaystyle\sum_{\emptyset\not=A\subset D}(-1)^{|A|-1}
\left(
\sum_{l=1}^{\infty} \sum_{k=-\infty}^{\infty} \displaystyle\bigvee_{i\in A}\left(a_{lki}\beta_{ji}^{1/\alpha_j}\right)
\right)^{\alpha_j}.$$

For the bivariate tail dependence it holds

$$ \lambda_{\{s,t\}}^{\left(C_{{\bf Y}}\right)}=2-\sum_{l=1}^{\infty} \sum_{k=-\infty}^{\infty}\displaystyle\sum_{j=1}^q\left(a_{lks}\beta_{js}^{1/\alpha_j}\vee a_{lkt}\beta_{jt}^{1/\alpha_j}\right)^{\alpha_j},$$

which, for the case $q=1$ leads to the result $ \lambda_{\{s,t\}}^{\left(C_{{\bf Y}}\right)}=2-\displaystyle\sum_{l=1}^{\infty} \displaystyle\sum_{k=-\infty}^{\infty}\left(a_{lks}\vee a_{lkt}\right)$ in Heffernan et al. (2007, \cite{hef+}).

%%%%%%%%%%%%%%%%%%%%%%%%%%%%%%%%%%%%%%%%%%%%%%%%%%%%%%%
\biblio{

\bibitem{cuadras+auge} Cuadras, C.M. and Aug\'e, J (1981). A continuous general multivariate distribution and its properties. Comm. Statist. A - Theory Methods 10, 339-353.

\bibitem{fougeres+nolan+rootzen} 
Foug\`eres, A.-L., Nolan, J.P. and Rootz\'en, H. (2009). Models for dependent extremes using scale mixtures. Scandinavian Journal of Statistics 36, 42-59.

\bibitem{hef+} Heffernan, J. E., Tawn, J. A., Zhang, Z. (2007). Asymptotically
(in)dependent multivariate maxima of moving maxima processes,
Extremes, 10, 57-82.

\bibitem{joe} Joe, H. (1997). Multivariate Models and Dependence Concepts.
Chapman \& Hall, London.

\bibitem{joe+hu} 
Joe, H. and Hu, T. (1996). Multivariate distributions from mixtures of max-infinitely divisible distributions. Journal of Multivariate Analysis 57, 240-265. 

\bibitem{li} 
Li, H. (2009). Orthant tail dependence of multivariate extreme value distributions.  Journal of Multivariate Analysis 100, 243-256. 

\bibitem{marshall+olkin} Marshall, A.W., Olkin, I. (1988). Families of multivariate distributions, J. Amer. Statist. Assoc. 83,
 834-841.
 
\bibitem{nelsen} Nelsen R. B. (1999). An Introduction to Copulas. Springer, New York.

\bibitem{schmid/schmidt} Schmid, F., Schmidt, R. (2007). Multivariate conditional versions of Spearman's rho and related measures of tail dependence. J. Multivariate Anal., 98, 1123-1140.

\bibitem{smith+weissman} Smith, R.L., Weissman, I. (1996). Characterization and estimation of the
multivariate extremal index. Technical Report, Univ. North,
Carolina.

\bibitem{tawn} 
Tawn, J. (1990). Modelling multivariate extreme value distributions. Biometrika 77, 2, 245-253.
}

\end{document}